\documentclass{article}

\usepackage{times}
\usepackage{epsfig}
\usepackage{graphicx}
\usepackage{amsmath}
\usepackage{amssymb}
\usepackage{algorithm}
\usepackage{amsthm}
\usepackage{amsfonts}
\usepackage{algorithmic}
\usepackage{color}
\newtheorem{thm}{Theorem}

\newtheorem{cor}{Corollary}

\usepackage{hyperref}
\usepackage{geometry}

\begin{document}
\title{Convex Shape Representation with Binary Labels for Image Segmentation: Models and Fast Algorithms}
\author{Shousheng Luo\\
School of Mathematics and Statistics, Henan University\\
Kaifeng, China.\\
{\tt\small sluo@henu.edu.cn}
\and
Xue-Cheng Tai\\
Department of Mathematics, Hong Kong Baptist University, \\ Kowloon Tong, Hong Kong\\
{\tt\small xuechengtai@hkbu.edu.hk}
\and
Yang Wang \\
Department of Mathematics, Hong Kong University of Science and Technology\\
Clear Water Bay, Kowloon,
Hong Kong
\\
{\tt\small yangwang@ust.hk}
}

\date{}
\maketitle
\begin{abstract}
   We present a novel and effective binary representation  for convex shapes. We show the equivalence between the shape convexity and some properties of the associated indicator function. The proposed method has two advantages. Firstly, the representation is based on a simple inequality constraint on the binary function rather than the definition of convex shapes, which allows us to obtain efficient algorithms for various applications with convexity prior. Secondly, this method is independent of the dimension of the concerned shape. In order to show the effectiveness of the  proposed representation approach, we incorporate it with a probability based model for object segmentation with convexity prior. Efficient algorithms are given to solve the proposed models using Lagrange multiplier methods and linear approximations. Various experiments are given to show the superiority of the proposed methods.
\end{abstract}
\section{Introduction}
Image segmentation with shape priors has attracted
much attention recently.
It is well-known that image segmentation plays a very important role in many modern applications. However, it is a very challenging task to segment the objects of interest accurately and correctly for
low quality images suffered from heavy noise, illumination bias, occlusions, etc.
Therefore, various shape priors are incorporated to improve the segmentation accuracy.

There is a long history of image segmentation with shape priors.
Early investigations on this topic are about object segmentation with concrete shape priors, see
\cite{Chan2005Level,Cremers2003Towards,Leventon2003Statistical}.
Nowadays, one pays more and more attentions on generic shape priors, such as connectivity \cite{vicente2008graph}, star shape \cite{gulshan2010geodesic,veksler2008star} and convexity \cite{Gorelick2017Convexity}.
More recently, an interesting research to segment an object consisting of several parts with given shape prior
drew a lot of attentions \cite{Zafari2016,isack2018k,isack2016hedgehog,Mahabadi2015Segment,toranzos2004sets}.
In this paper, we focus on the topic of image segmentation with convexity prior.


{\bf Related works}~
Star shape is closely related to convexity, and is one of the widely investigated priors in the literature.
A region is called star shape with respect to a given point if all the line segments between all points in this region and the
referred point belong to this region as well. As far as we know,
it was first adopted as shape prior for image segmentation in \cite{veksler2008star}.
Then it was extended to star shape with more than one referred points and geodesic star shape \cite{gulshan2010geodesic}.
In \cite{yuan2012fast,yuan2012efficient}, efficient algorithm based on graph-cut was proposed.
Recently, star shape prior was encoded in neural network for skin lesion segmentation \cite{mirikharaji2018star}.

Recently, image segmentation with convexity prior attracted increasing attentions for low quality images,
and various methods were proposed in the literature. These methods can be categorized into two classes.
The first one is based on level set function \cite{Bae2017Augmented,Ukwatta2013Efficient}, and the other is based on binary representation \cite{gorelick2017multi,Gorelick2017Convexity}.

Convexity prior was investigated via level set method in \cite{Ukwatta2013Efficient}. Then this idea was adopted in \cite{Bae2017Augmented,yang2017level}. Actually, these methods utilize the fact that the curvature of the convex shape boundary must be nonnegative. More recently, this idea was developed in \cite{Luo:2018tl,Yan2018}, where the authors extended
the nonnegative curvature on the boundary (zero level set curve) to all the level set curves in the image domain.
In addition, efficient algorithms were also proposed thanks to the fact that the curvature on all level set curves
can be computed by the Laplacian of the associated signed distance function.
This method was extended to multiple convex objects segmentation in \cite{luo2019Cnvx} and applied for convex hull problems in  \cite{li2019variational}.

Binary representation has also been  used for image segmentation with convexity prior. These methods utilize the definition of convex regions.
In \cite{Gorelick2017Convexity}, 1-0-1 configurations on all lines in the image domain were  penalized
to promote convex segmentation, where 1 (resp. 0) is used to represent that the corresponding point belongs to foreground (resp. background).
Then this method was extended for multiple convex objects segmentation in \cite{gorelick2017multi}.
Algorithms based on trust region and graph-cut algorithms were studied in \cite{Gorelick2017Convexity} and \cite{gorelick2017multi} using linear or quadratic approximations, respectively.
For a convex region, 
the line segments between any two points should not pass through the object boundary, which is used to characterize convexity \cite{royer2016convexity}.
This description was incorporated into the multicut problem for image segmentation \cite{royer2016convexity}, and an algorithm based on branch-and-cut method was used to solve this problem.

The two approaches mentioned above suffer from some disadvantages.
For the level set approach, it is not easy  to extend the method for three dimensional (3D) image segmentation with convexity prior.
Although one can use the nonnegativity of Gaussian curvature to characterized the convexity of 3D objects \cite{elsey2009analogue},
it is a very challenging task to solve the corresponding models with very complex constraints involving curvatures. For binary approaches,
the computational cost is also an insurmountable problem for the method in \cite{Gorelick2017Convexity}.
Therefore, it is  necessary  to develop new methods for convex shape representation, which
can be extended for 3D convex object representations.

{\bf Contributions}~
In this paper, we propose a novel binary label method for convex shape representation.
Let us consider a convex shape in $\mathbb{R}^d$.
For any given ball in $\mathbb{R}^d$ (disc for $d=2$) centered on the boundary,
the volume (area for $d=2$) of the ball inside (resp. outside) the convex region is less (resp. greater)
than half of the ball volume (see figure \ref{fig:Interp}).
Obviously, the conclusion is also true for balls centered outside the object region.
Therefore, we obtain an easy and simple equivalence between convex shapes and their binary representations, which is regardless of the dimension of the objects.

According to the observation above, we can develop an efficient binary representation for convex objects,
which is  an inequality constraint on the indicator function associated to the concerned shape.
Let $b$ be a positive radial function defined on a given ball with integral equalling $1$, and
$u$ be the associated indicator function with the considered object, i.e. $u=1$ outside the object and $0$ inside the object.
Finally, the object
 is convex if and only if  $ub\ast u\geq0.5u$ for all radial functions $b$ (the derivation details will be given in Section \ref{sec:BR}), where $\ast$ denotes convolution operator in $\mathbb{R}^d$.
Accordingly, we obtain an equivalent description for convex shape based on binary representation by imposing an inequality constraint. In addition, this method can be easily extended to multiple convex objects representation using the technique in
\cite{luo2019Cnvx}.

Comparing to the methods in \cite{Gorelick2017Convexity,royer2016convexity},
the advantages of the proposed method can be summarized as follows:
\begin{enumerate}
    \item
The proposed method is a very general convex shape representation technique, which is regardless of the object dimension.
\item
Simplicity and numerical efficiency is another advantage of our method.
The proposed method is very simple, which allows us to design efficient algorithms. In this work, we use these techniques for image segmentation. It is easy to extend these ideas for other applications with general shape optimization problems with convex shape prior.
\end{enumerate}


In order to show the effectiveness of the proposed method, we apply
it to image segmentation with convexity prior,
although it can be used for general problems with convexity prior, e.g. convex hull \cite{li2019variational}.
In this paper, the proposed convexity representation method is incorporated into a probability based model for image segmentation.
The region force term is computed as the negative log-likelihood of probabilities belonging to the foreground and background, where the probabilities are fitted by mixed Gaussian method \cite{Rother:2004}.

An efficient algorithm for the proposed model is developed using Lagrange multiplier method.
Firstly, we can write down the associated Lagrange function of the
segmentation model with the inequality constraint for convexity.
Secondly, we use the technique in \cite{liu2011fast,wang2017efficient,wang2009edge} to approximate the boundary length or area regularization.
For the proposed iterative algorithm, explicit binary solution is available for each step after linearizing the quadratic constraint. As for the multiplier update, it is updated by gradient ascent method which is simple and also turns out to be efficient.

In order to improve the stability of the algorithm, more techniques are added to our algorithm.
\romannumeral1) The binary function is updated only on a narrow band of the boundary of the current object estimate.
\romannumeral2) The value of the binary function on the boundary is set to 0.5 for the area computation in the implementation because
the measure of boundary is zero in continuous setting, but not zero in discrete setting.

The rest of the paper is structured as follows. We will present the binary representation method for convex shape and give details of the image segmentation model with convexity prior in Section \ref{sec:BR}.
Numerical algorithms are proposed in Section \ref{sec:Num}.
Some experimental results are  demonstrated in Section \ref{sec:Exp}.
We conclude this paper and discuss future works in Section \ref{sec:Con}.

\section{The proposed method}\label{sec:BR}
In this section, we will present the proposed binary representation for convex objects, and then
incorporate it with probability models for image segmentation.
\subsection{Binary representation of convex object}
Before presenting the binary representation for convexity shapes, we introduce some
notations firstly.
For a given set $S\subset\mathbb{R}^d(d\geq2)$,
$S^c$ denotes the complementary set of $S$, and
the associated indicator function $v$ with $S$ is defined as
\begin{equation}
    v(x)=\left\{\begin{array}{ll}
    1,x\in S^c,\\
    0,x\in S.
    \end{array}
    \right.
\end{equation}
Let $B_r(x)\subset\mathbb{R}^d$ denote the ball centered at $x$ with radius $r>0$, and $b_r$ denotes the
radial function with $\int_{B_r(0)}b_r(x)dx=1$ and $b_r(x)\geq0$ if $|x|\leq r$ and $b_r(x)=0$ if $|x|>r$.

\begin{figure}
    \centering
    \includegraphics[width=4cm,height=3cm]{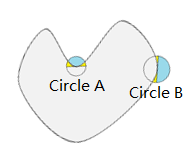}
    \caption{Interpretation for Theorem \ref{thm1}. }
    \label{fig:Interp}
\end{figure}
\begin{thm}\label{thm1}
   Suppose $D\subset \mathbb{R}^d~(d\geq2)$ is the object region that we want to extract. 
    Then $D$ is convex if and only if the following inequality holds
   \begin{eqnarray}
     \mathcal{A}(D^c\bigcap B_r(x))\geq \frac{1}{2}\mathcal{A}(B_r(x))  \label{eq:char2}
   \end{eqnarray}
   for all $r>0$ and all $x~\text{on}~\partial D$,  where $\mathcal{A}$ denotes the measure (volume or area) in $\mathbb{R}^d$.
\end{thm}
The proof is very easy and will be presented in Appendix. An intuitive interpretation for two-dimensional case is illustrated in Figure \ref{fig:Interp}. Circle A (resp. B) is at a nonconvex (resp. convex) position on the boundary. Obviously, the
inequality in (\ref{eq:char2}) holds for $x$ on the convex points, and it is violated on  nonconvex points.
Inequality (\ref{eq:char2}) provides a boundary-based characterization for convex region.
In fact, the inequality is true for all  $x\in D^c$, i.e.
\begin{equation}
    \mathcal{A}(D^c\bigcap B_r(x))\geq \frac{1}{2}\mathcal{A}(B_r(x)), x \in D^c \label{eq:char3}.
\end{equation}
%

Let $u$ be the indicator function of $D$.
According to the notations above, inequality (\ref{eq:char3}) is equivalent to
\begin{equation}
    \int_{D^c}\!b_r(y-x)dy\!=\!\int_{\Omega}\!u(y)b_r(y-x)dy\geq \frac{1}{2},
     \ x\in D^c.
    \label{eq:char4}
\end{equation}
 Based on the discussions above, we can obtain the following equivalent description for the convexity of $D$.
\begin{cor}\label{cor1}
Under the assumption in Theorem \ref{thm1},
$D$ is convex if and only if
\begin{equation}
   \mathcal{C}_r(u)= u(x)b_r\ast u(x)-\frac{1}{2}u(x)\geq 0 \label{eq:char5}
\end{equation}
for all $x\in \Omega$, where $b\ast u(x)=\int_{\Omega}u(y)b_r(y-x)dy$ denotes the convolution in $\mathbb{R}^d$.
\end{cor}
The proofs of (\ref{eq:char4}) and Corollary \ref{cor1} will be presented in the appendix.

\subsection{Image segmentation model with convex prior}\label{sec:segC}
In this section we will incorporate the proposed binary representation with a
segmentation model for convex object segmentation.
Let $I$ be a given image $I: \Omega \mapsto \mathbb{R}^p$ defined on $\Omega\subset \mathbb{R}^d$ ($p=1$ for gray image and $3$ for color image), and $D\subset \Omega$ be the object domain of interest to extract.
The Potts segmentation model with boundary length regularizer and convexity prior can be written as
\begin{equation}
\begin{aligned}
   \min_D\int_{D}w_0f_{0}(x)dx\!+\!\int_{D^c} w_1f_{1}(x)dx+\lambda|\partial D|,\\
   \hfill\text{subject to }~~D~~\text{being convex},
\end{aligned}\label{eq:genalM}
\end{equation}
where $f_{0}$ (resp. $f_1$) are some given similarity measures on the object $D$ (resp. background $D^c$) and $w_0,w_1,\lambda>0$ are user-specified trade-off parameters.

Let $u$ be the indicator function of $D$. As pointed in \cite{mirpalpar2007}
the length or area of boundary $D$ can be approximated by
\begin{equation}
    |\partial D|\approx \mathcal{L}_\sigma (u)=\sqrt{\frac{\pi}{\sigma}}\int_{\Omega}u(x)G_\sigma\ast{(1-u)}(x)dx,
\end{equation}
when $0<\sigma\ll 1$, where $G_\sigma$ is the Gaussian kernel
\begin{equation}
    G_{\sigma}(x)=\frac{1}{(4\pi\sigma)^{d/2}}\exp\left(-\frac{|x|^2}{4\sigma}\right).
\end{equation}

According to the results in last section, the Potts model (\ref{eq:genalM}) can be equivalently formulated as
\begin{equation}
\begin{aligned}
    \min_{u\in\{0,1\}}\!\int_{\Omega}[w_1\!f_{1}-w_0\!f_{0}]udx\!+\!\lambda\mathcal{L}_\sigma(u),\mathcal{C}_r(u)\geq0,\label{eq:ModBny1}
\end{aligned}
\end{equation}
for all $r>0$. One can finds more details in \cite{esedog2015threshold,wang2009edge}.

In fact, it is not necessary to impose the constraint $\mathcal{C}_r(u)\geq0$ for all $r>0$.
In our implementations, the segmentation model with convexity prior (\ref{eq:ModBny1}) is simplified as
\begin{equation}
\begin{aligned}
    \min_{u\in\{0,1\}}\!\int_{\Omega}fudx\!+\!\lambda\mathcal{L}_\sigma(u),\mathcal{C}_{r_i}(u)\geq0,\label{eq:ModBny}
\end{aligned}
\end{equation}
for some given  $r_i>0  \ (i=1,2,\cdots,n)$, where $f=w_1f_1-w_0f_0$ and $n$ is a user-specified integer.

\subsection{Labels and region force}
For complex images, we need to use semi-supervised segmentation model.
We assume the labels of some points from the foreground and background already  have been known.
Let $R_{ob}$ and $R_{bg}$ be the sets of points with known labels. The semi-supervised  segmentation model (\ref{eq:ModBny}) can be written as
\begin{equation}
\begin{aligned}
    \min_{u\in\{0,1\}}\!\int_{\Omega}fudx+\lambda\mathcal{L}_\sigma(u), \mathcal{C}_r(u)\geq0,    u\in L, \label{eq:ModBnyL}
\end{aligned}
\end{equation}
where $L=\{u|u(x)=1,x\in R_{bg},u(x)=0,x\in R_{ob}\}$.

The region force $f$ (or $f_0,f_1$) plays a very important role.
Here we adopt the probability method to compute $f_0$ and $f_1$.
For each $x\in\Omega$ and the given image $I$, assume $p_0(x|I(x)),p_1(x|I(x))$ are the estimated probabilities for the foreground and background.  Then we use the negative log-likelihood of $p_i$ as $f_i$ ($i=0,1$), i.e.
\begin{equation}
    f_{i}(x)=-\ln(p_{i}(x|I(x))). \label{eq:Frc}
\end{equation}
As for the probabilities $p_0$ and $p_1$, we use mixed Gaussian method in \cite{Rother:2004} to
estimate them. Given estimates of foreground and background,
the density functions of $I(x)$ on foreground and background are fitted by
different mixed Gaussian distributions, i.e.
\begin{align}
    G_{0}(I(x))=\sum_{k=1}^{N_0}c_k^0G(I(x);\mu_k^0,\Sigma_k^0),\label{eq:p0}\\
    G_{1}(I(x))=\sum_{k=1}^{N_1}c_k^1G(I(x);\mu_k^1,\Sigma_k^1),\label{eq:p1}
\end{align}
where $c_k^i,\mu_k^i,\Sigma_k^i,k=1,2,\cdots,N_i, i=0,1$ are the fitting parameters for given $N_0$ and $N_1$.
Here $c_k^i$ are the portions of different Gaussian distributions with mean $\mu_k^i$ and variance $\Sigma_k^i$.
We can estimate the probabilities for the foreground and background as
\begin{equation}
    p_0(x|I(x))=\frac{\gamma_0G_{0}(I(x))}{\gamma_0G_{0}(I(x))+\gamma_1G_{1}(I(x))}
\end{equation}
and $p_1(x|I(x))=1-p_1(x|I(x))$, where $\gamma_0$ and $\gamma_1$ are weighting parameters for the foreground and background probabilities.

\section{Numerical method}\label{sec:Num}
This section is devoted to numerical techniques for the proposed model (\ref{eq:ModBnyL}).
Hereafter, we denote $C_{r_i}(u)$ (resp. $b_{r_i}$) by $C_i$ (resp. $b_i$) for notational simplicity.
We can write down the Lagrange functional of (\ref{eq:ModBnyL}) as
\begin{equation}
\min_{u\in\{0,1\}}\max_{g_i\geq0}\int_{\Omega}\![fu\!-\!\sum_{i=1}^n g_i\mathcal{C}_i(u)]dx
\!+\!\lambda\mathcal{L}_\sigma(u),u\in L,\label{eq:UpdL1}
\end{equation}
where $g_i$ is the Lagrange multiplier associated with $C_i(u)\geq0$.


We can use alternating direction method to solve (\ref{eq:UpdL1}).
For given initialization $u^0$ and $g_i^0$, we use projection gradient ascent method to update $g_i$, i.e.
\begin{equation}
g_i^{t+1}=\text{Proj}^+(\tilde{g}_r^t)=\max\{0,\tilde{g}_r^t\}\label{eq:GmUpd}
\end{equation}
with $\tilde{\gamma}_i^t=g_i^t-\tau\mathcal{C}_i(u^t)$ for $i=1,2,\cdots,n$, where $\tau>0$ is the step size.

For the update of $u$,
we linearize the last two terms at $u^k$ in (\ref{eq:UpdL1}) to approximate the objective functional, i.e.
\begin{equation}
u^{t+1}=\arg\min_{u\in\{0,1\}}\!\int_{\Omega}F(x;u^t,g_i^t)udx,u\in L,\label{eq:UpdL2}
\end{equation}
where
\begin{align}
    F^k(u^t,g_i^t)=\sum_{i=1}^n[0.5g_i^k-b_i\ast(u^t+g^t_iu^t)]\notag\\
    +f+\lambda G_{\sigma}\ast(1-2u^t) .
    \label{eq:UpdApp}
\end{align}
Therefore, the solution  $u^{t+1}$ is given by the following explicit formula:
\begin{equation}
    u^{t+1}(x)=\text{Proj}_L(\tilde{u}^t)=
    \left\{\begin{array}{ll}
    1& x\in R_{bg},\\
    0&x\in R_{ob},\\
    \tilde{u}^{t}& \text{otherwise},
    \end{array}
    \right. \label{eq:updF}
\end{equation}
where
\begin{equation}
    \tilde{u}^{t}(x)=\left\{\begin{array}{ll}
    1& F(x;u^t,g_i^t)\leq0,\\
    0&F(x;u^t,g_i^t)>0.\\
    \end{array}
    \right.\label{eq:utld}
\end{equation}

Finally, the algorithm for (\ref{eq:UpdL1}) is summarized as Algorithm \ref{Alg:U}.
\begin{algorithm}
\caption{}\label{Alg:U}
\begin{algorithmic}
\STATE Initialize: $u^0, g_i^0$, $t=0$, and maximum number $T>0$
\STATE While $t<T$ and termination criterion on $u$ false
\STATE ~ Update $u^{t+1}$ using (\ref{eq:updF});
\STATE ~ Update $g_i^{t+1}$ using (\ref{eq:GmUpd});
\STATE ~ $t=t+1$;
\STATE End(while)
\end{algorithmic}
\end{algorithm}

\subsection{Initialization}
We adopt the following method to initialize $u$ and region force $f$.
We use the convex hull  of the subscribed labels $R_{ob}$, denoted by $D=CH(R_{ob})$, as
the initial foreground, 
and the indicator function of $D=CH(R_{ob})$ as the initialization of $u$.

As for the estimate of the region force term, we use the following method to fit the distributions $G_0$ and $G_1$ in (\ref{eq:p0}) and
(\ref{eq:p1}), and compute the region force using (\ref{eq:Frc}).
According to initialization of foreground, the distribution $G_0$ is fitted using $I(x)$ for $x\in D=CH(R_{ob})$.
In order to obtain a more accurate initialization of $G_1$, we use $I(x)$ on $D^{c,s}$ to estimate it, where
\begin{equation}
D^{c,s}=\{x\in \Omega|\min_{y\in D}\|x-y\|_2>s\},\label{eq:Dc}
\end{equation}
where $s>0$ is a parameter.


\subsection{Numerical details}
For a given $M\times N$ digital  image $I$, we just view it as a discrete image  defined on $\Omega=[0,M-1]\times [0,N-1]$ with mesh size
$h=1$.
One detail deserving our attention is the convolution operation $b_r\ast u$ in discrete implementation.
Let $D^t=\{x|u^t(x)=0\}$ be the estimated foreground  corresponding to the current binary function $u^t$.
In order to compute the convolution $b_r\ast u^t$ accurately, we set $u^{t}=0.5$ on $\partial D^{t}$ , which
is extracted by the matlab function \emph{bwperim} in the implementation.

In order to improve the stability of the algorithm, several techniques will be used.
Firstly, we replace the constraint $u(x)=0$ for all $x\in R_{ob}$  to
$u(x)=0$ for all $x\in CH(R_{ob})$, and require $u(x)=1$ for $x\in R_{bg}\backslash CH(R_{ob})$ if
$R_{bg}\bigcap CH(R_{ob})\neq\emptyset$.
Secondly, the update of $u^{t+1}$ is constrained on a narrow band of the current estimated boundary
\begin{equation}
    S(u^t)=\{x||b_{r_0}\ast u^t-u^t|\geq \rho\},\label{eq:Nrw}
\end{equation}
where $\rho>0$ is an user-specified parameter.
In addition, a proximal term $\frac{\theta}{2}\|u-u^t\|$ is added into the objective functional for the update of $u$.

Using the binary constraint on $u^{t+1}$,
(\ref{eq:utld}) is improved as
\begin{equation}
    \tilde{u}^{t}(x)=\left\{\begin{array}{ll}
    1& F(x;u^t,g_i^t)+\theta(0.5-u^t)\leq0,\\
    0&F(x;u^t,g_i^t)+\theta(0.5-u^t)>0.\\
    \end{array}
    \right.\label{eq:utld1}
\end{equation}
Therefore, the update for $u$ can be modified as
\begin{align}
    u^{t+1}(x)=\left\{\begin{array}{ll}
    \text{Proj}_L(\tilde{u}^t)(x)& x\in S(u^t),\\
    u^t(x)&\text{otherwise},
    \end{array}
    \right.\label{eq:UpdF}
\end{align}
where $\tilde{u}^t$ is as in (\ref{eq:utld}) and
\begin{equation}
   \text{Proj}_L(\tilde{u}^t)=
    \left\{\begin{array}{ll}
    1& x\in R_{bg}\backslash CH(R_{ob}),\\
    0&x\in CH(R_{ob}),\\
    \tilde{u}^{t}& \text{otherwise}.
    \end{array}
    \right. \label{eq:updF1}
\end{equation}


In addition,
we terminate the iteration  when
the relative variation between two iterations is less than a tolerance.  Let
$v_1,v_2$ be two binary functions. The relative variation of $v_2$ with respect to $v_1$ is
defined as by
\begin{equation}
    R(v_1,v_2)=\frac{\int_{\Omega}|v_2-v_1|dx}{\int_{\Omega}v_1dx}.
\end{equation}
In order to avoid early termination inappropriately, we compute the relative variation every 300 iterations, i.e. $R(u^{t},u^{t-300})$.
In addition, the region force term is updated in the iterative procedure only when it is needed.
In this paper, we update $f$ using the current estimate foreground $\{x|u^t=0\}$ and background $\{x|u^t(x)=1\}$ every $50$ iterations.

According to implementation details above, we can summarize the algorithm for the proposed method as Algorithm \ref{Alg:detail}.
\begin{algorithm}
\caption{Algorithm for the proposed model}\label{Alg:detail}
\begin{algorithmic}
\STATE1. Input: The concerned image and subscribed labels $R_{bg}$ and $R_{ob}$;
\STATE2. Initialization: $u^0$ and the region force term $f$;
\STATE3. Set: $Rv>\epsilon$, $t=0$, maximum iteration number $T>0$ and tolerance $\epsilon>0$, and
    $g_{i}^0\in\mathbb{R}^{M\times N}$, $i=1,2,\cdots,n$.
\STATE4. While $t<T$\& $Rv>\epsilon$
\STATE5. ~~Determine narrow band $S(u^t)$ by (\ref{eq:Nrw}) ;
\STATE6. ~~Update $u^{t+1}$ using (\ref{eq:UpdF});
\STATE7. ~~Update $g_{i}^{t+1}$ using (\ref{eq:GmUpd}) for $i=1,2,\cdots,n$;
\STATE8. ~~$t=t+1$;
\STATE9. ~~Update $f$ using current estimate if $mod(t,50)=0$;
\STATE10. ~Compute  $Rv=(u^t,u^{t-300})$ if $mod(t,300)=0$;
\STATE11. End(while)
\end{algorithmic}
\end{algorithm}
%

We can see that the main operations in Algorithm \ref{Alg:detail} are
convolutions,  which are $b_i\ast u$ for convexity constraint and narrow band determination and $G_\sigma\ast u^t$ for boundary measure approximation. It is well-known that
the convolution operation can be implemented by FFT efficiently.
Therefore, the proposed algorithm is very cheap and efficient.

\section{Experiments}\label{sec:Exp}
In this section, we will present some numerical results to show the efficiency and effectiveness of the proposed method and algorithm.
A lot of experiments were conducted on various images, and the results show the effectiveness of the proposed method in
preserving the convexity of shapes. Here we only demonstrates some of them. 

In the implementation, some parameters are kept the same for all the experiments for simplicity.
We set $w_0=w_1=0.5$ and $\lambda=0.1$ in the model (\ref{eq:ModBnyL}).
As for the boundary length approximation term, we use $5\times5$ Gaussian kernel with variance 0.5 generated by matlab bulit-in function $fspecial$.
The integer $N_0$ (resp. $N_1$) of the Gaussian distributions for foreground (resp. background) is set to $2$ (resp. $3$), and $\gamma_0,\gamma_1$ are set to $1$.
The parameter $s$ in (\ref{eq:Dc}) and $\rho$ in (\ref{eq:Nrw}) are set to $5$ and $2$, respectively.
The proximal parameter $\theta$ in (\ref{eq:utld1}) equals to $1$ in the implementation.
The step size for dual variable $g_i$ update is set to $1$.
In numerical implementation we choose $n=4$ and  $r_i=4+5(s-1)$ for $s=1,2,3,4$ for the convexity constraint $C_i(u)\geq0$ if they are not specified.
As for radial functions, we just choose $b_{r_i}$ as the uniform function on the discs $B_{r_i}$, i.e. $b_r(x)=\frac{1}{\pi r^2}$, which is
generated by the matlab built-in function \emph{fspecial}.
The radius $r_0$ for the narrow band determination is set to $3$.
The relative variation tolerance $\epsilon$ and the maximum iteration number are set to $0.001$ and $5000$.
\begin{figure}
    \centering
    \includegraphics[width=2.4cm]{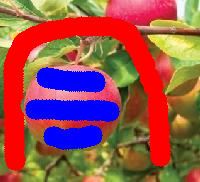}
    \includegraphics[width=2.4cm]{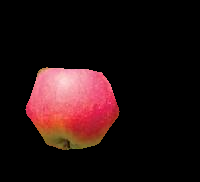}
    \includegraphics[width=2.4cm]{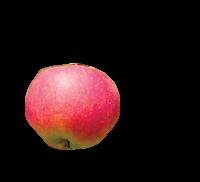}\\
    \includegraphics[width=2.4cm]{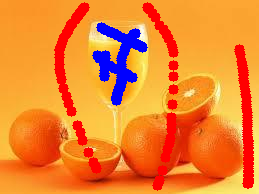}
    \includegraphics[width=2.4cm]{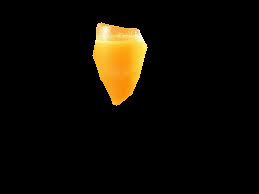}
    \includegraphics[width=2.4cm]{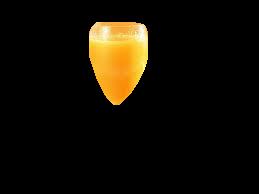}\\
    \includegraphics[width=2.4cm,height=2cm]{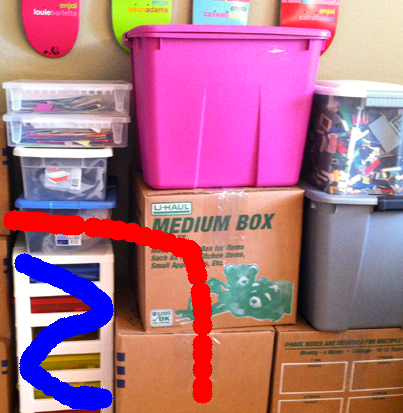}
    \includegraphics[width=2.4cm,height=2cm]{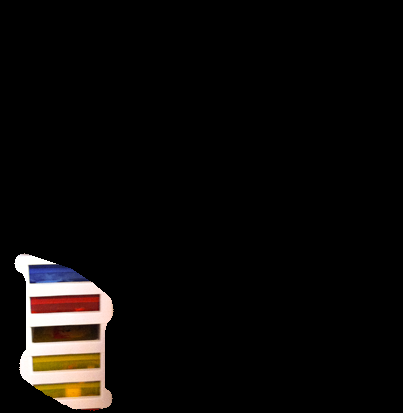}
    \includegraphics[width=2.4cm,height=2cm]{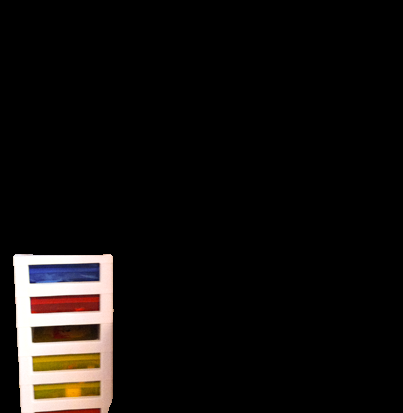}\\    \includegraphics[width=2.4cm]{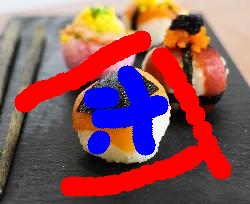}
    \includegraphics[width=2.4cm]{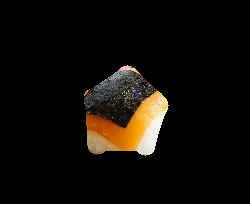}
    \includegraphics[width=2.4cm]{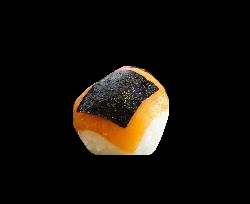}\\
    \includegraphics[width=2.4cm]{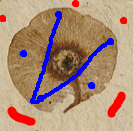}
    \includegraphics[width=2.4cm]{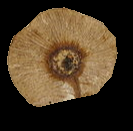}
    \includegraphics[width=2.4cm]{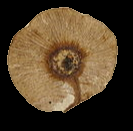}\\
    \includegraphics[width=2.4cm]{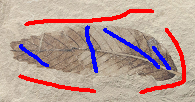}
    \includegraphics[width=2.4cm]{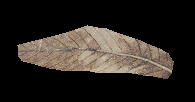}
    \includegraphics[width=2.4cm]{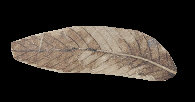}\\

    \includegraphics[width=2.4cm,height=2cm]{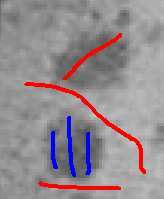}
    \includegraphics[width=2.4cm,height=2cm]{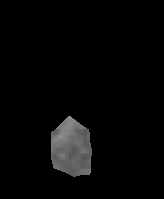}
    \includegraphics[width=2.4cm,height=2cm]{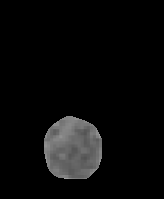}\\
    \includegraphics[width=2.4cm]{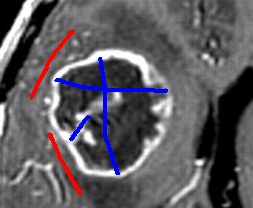}
    \includegraphics[width=2.4cm]{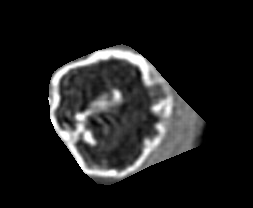}
    \includegraphics[width=2.4cm]{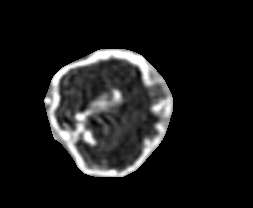}\\
    \includegraphics[width=2.4cm,height=2cm]{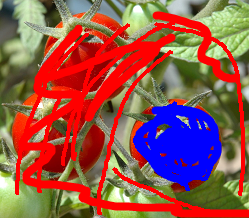}
    \includegraphics[width=2.4cm,height=2cm]{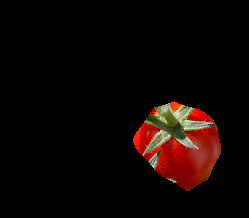}
    \includegraphics[width=2.4cm,height=2cm]{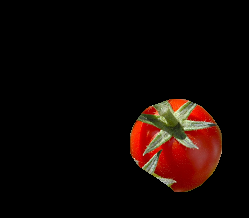}\\
    \includegraphics[width=2.4cm,height=2cm]{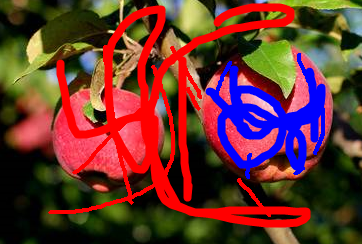}
    \includegraphics[width=2.4cm,height=2cm]{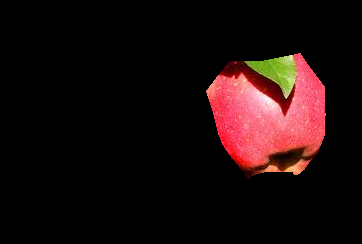}
    \includegraphics[width=2.4cm,height=2cm]{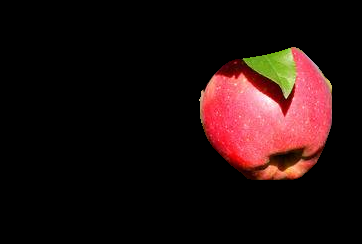}
\caption{Comparison: Input images with pre-labled pixles(left), results by the method in \cite{Gorelick2017Convexity} (middle) and by the proposed method (right).}\label{fig:compob1}
\end{figure}
\begin{figure}
\centering
\includegraphics[width=2.5cm]{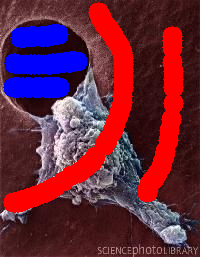}
\includegraphics[width=2.5cm]{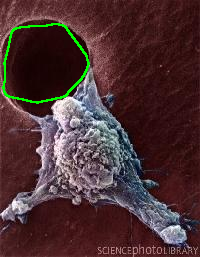}
\includegraphics[width=2.5cm]{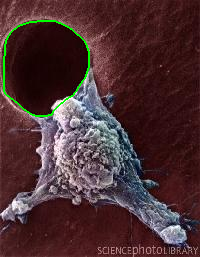}\\
\includegraphics[width=2.5cm,height=2cm]{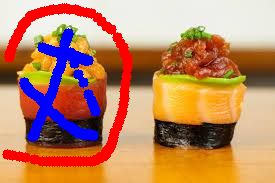}
\includegraphics[width=2.5cm,height=2cm]{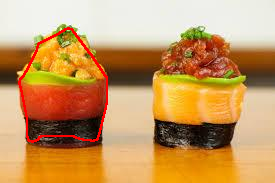}
\includegraphics[width=2.5cm,height=2cm]{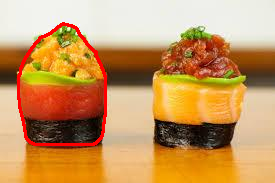}
\\
\includegraphics[width=2.5cm,height=1.8cm]{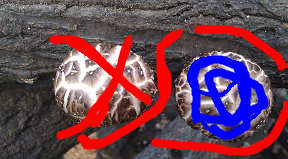}
\includegraphics[width=2.5cm,height=1.8cm]{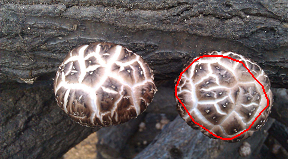}
\includegraphics[width=2.5cm,height=1.8cm]{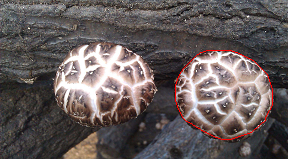}\\
\includegraphics[width=2.5cm]{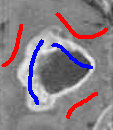}
\includegraphics[width=2.5cm]{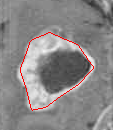}
\includegraphics[width=2.5cm]{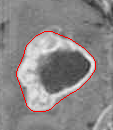}
\caption{Comparison: Input images with pre-labled pixles(left), results by the method in \cite{Gorelick2017Convexity} (middle) and by the proposed method (right). 
}\label{fig:compEg1}
\end{figure}
\begin{figure}
\centering
    \includegraphics[width=2.5cm]{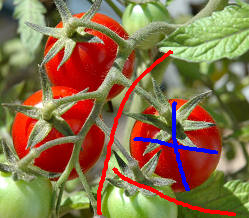}
    \includegraphics[width=2.5cm]{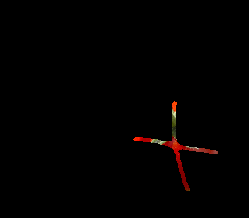}
    \includegraphics[width=2.5cm]{img/img27/segObj.png}\\
    \includegraphics[width=2.5cm]{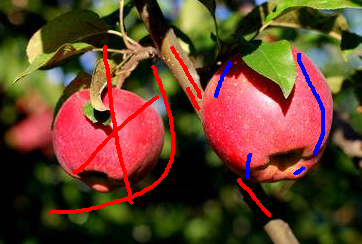}
    \includegraphics[width=2.5cm]{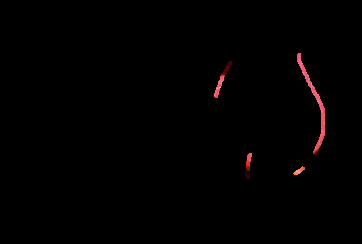}
    \includegraphics[width=2.5cm]{img/img33/segObj.png}\\
    \includegraphics[width=2.5cm]{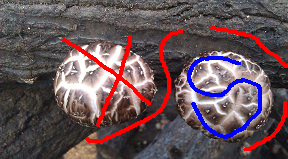}
    \includegraphics[width=2.5cm]{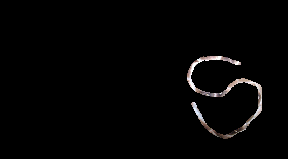}
    \includegraphics[width=2.5cm]{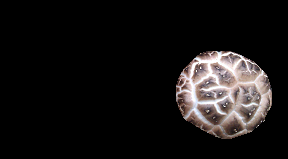}
\caption{Comparison: results with fewer pre-given labels.}\label{fig:compL}
\end{figure}

\subsection{Result comparison}
Some experimental results are presented in figures \ref{fig:compob1}, \ref{fig:compEg1} and \ref{fig:compL} to compare the proposed method with the method in \cite{Gorelick2017Convexity}.
For the method \cite{Gorelick2017Convexity}, we use $11\times 11$ stencil  and penalty parameter equaling to $2$.

Some of the testing images and results of the method \cite{Gorelick2017Convexity} are downloaded from website \url{http://vision.csd.uwo.ca/code/}.
In order to compare the results, the segmentation objects are extracted (see figure \ref{fig:compob1}) or the segmentation object boundaries are drawn (see figure \ref{fig:compEg1}).
Besides the images (the top four (resp. two) images in figure \ref{fig:compob1} (resp. figure \ref{fig:compEg1}) ) downloaded from the website, some more experiments on other images conducted to compare the proposed method and  the method \cite{Gorelick2017Convexity}.


We can see that the proposed method is superior to the method \cite{Gorelick2017Convexity} by comparing the results in
figure \ref{fig:compob1}.
Although the method \cite{Gorelick2017Convexity} can extract the main parts of the concerned objects,
the proposed method can extract the objects more completely and accurately, e.g. the apple, the lotus leaf and the tomato.
Taking the third image in figure \ref{fig:compob1} as an example, we can see that two corners of the object are smeared by the method \cite{Gorelick2017Convexity}, while the result by the proposed method is more complete and accurate.
For the results in figure \ref{fig:compEg1}, we can see that the proposed method
can touch the object boundary precisely and accurately, while the method \cite{Gorelick2017Convexity} fails to
capture the object boundary accurately. For example, the result of the first image in figure \ref{fig:compEg1}, there are only
few points of the extracted boundary reach the concerned objects' boundary.

In addition, the method \cite{Gorelick2017Convexity} usually needs more pre-labeled pixels than the proposed method.
When we do not have enough pre-labeled pixels, the results by \cite{Gorelick2017Convexity}  are often very poor. The last two images in figure \ref{fig:compob1} are examples.
The method \cite{Gorelick2017Convexity} fails to get meaningful results with fewer labels (see figure \ref{fig:compL}).
Therefore, a lot of pre-given labels are needed  for the images (the tomato and apple images in figure \ref{fig:compob1} and the mushroom image in \ref{fig:compEg1}) to obtain
a meaningful segmentation for the method \cite{Gorelick2017Convexity}, although the proposed method does not need.

\subsection{Sensitivity to the radius}
\begin{figure}
    \centering
    \includegraphics[width=2.5cm]{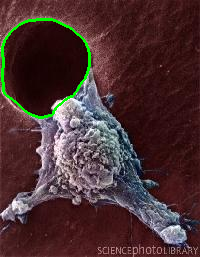}
    \includegraphics[width=2.5cm]{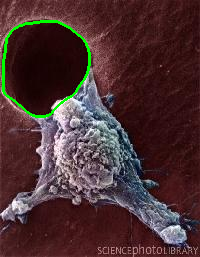}
    \includegraphics[width=2.5cm]{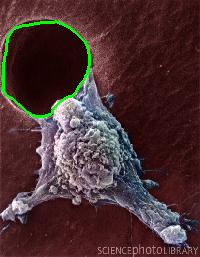}\\
    \includegraphics[width=2.5cm]{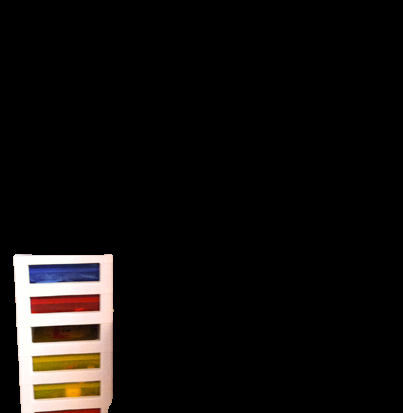}
    \includegraphics[width=2.5cm]{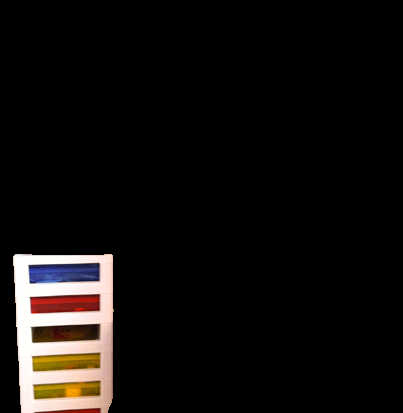}
    \includegraphics[width=2.5cm]{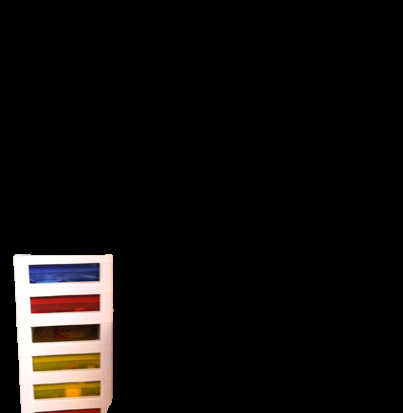}\\
    \includegraphics[width=2.5cm]{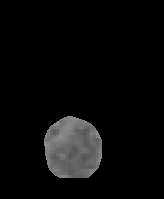}
    \includegraphics[width=2.5cm]{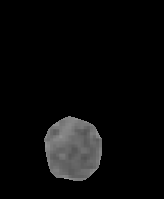}
    \includegraphics[width=2.5cm]{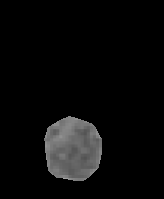}\\
    \includegraphics[width=2.5cm]{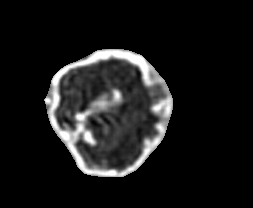}
    \includegraphics[width=2.5cm]{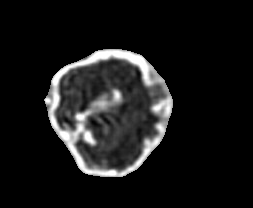}
    \includegraphics[width=2.5cm]{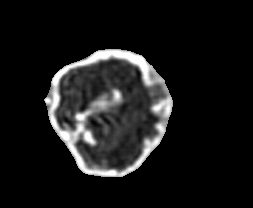}
   \\
    \includegraphics[width=2.5cm]{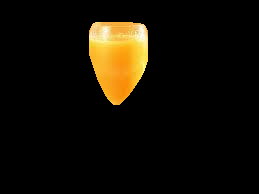}
    \includegraphics[width=2.5cm]{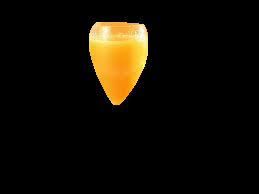}
    \includegraphics[width=2.5cm]{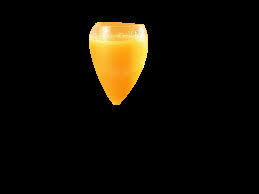}
  \\
    \includegraphics[width=2.5cm]{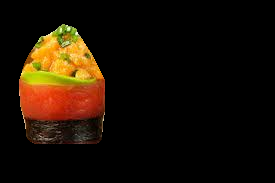}
    \includegraphics[width=2.5cm]{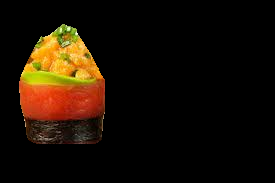}
    \includegraphics[width=2.5cm]{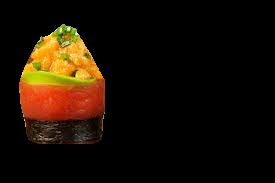}\\
    \includegraphics[width=2.5cm]{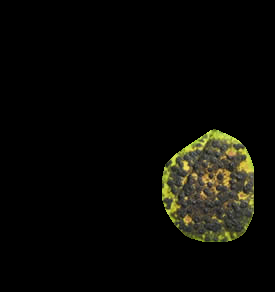}
    \includegraphics[width=2.5cm]{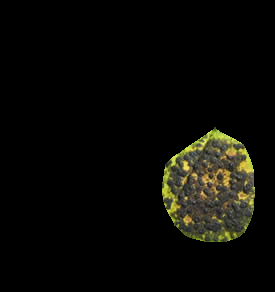}
    \includegraphics[width=2.5cm]{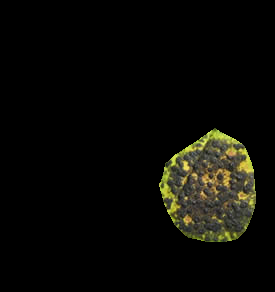}\\
    \includegraphics[width=2.5cm]{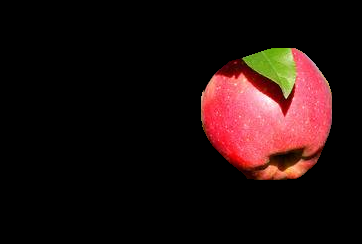}
    \includegraphics[width=2.5cm]{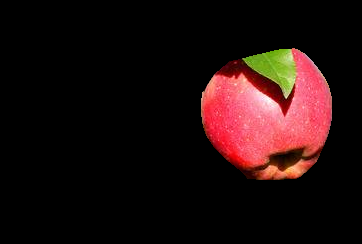}
    \includegraphics[width=2.5cm]{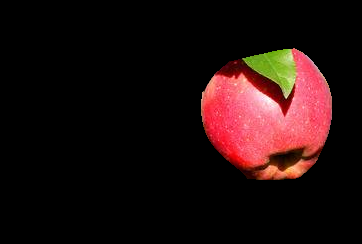}
    \caption{Sensitivity to radius of the radial functions.}\label{fig:rad}
\end{figure}

Some experiments with  different radial functions were conducted to investigate the robustness of the proposed method. The results are presented in figure \ref{fig:rad}. 
The radius of the radial functions for the images from left to right in
figure \ref{fig:rad} are $[4,10,14,20]$, $[5,10,15,24]$ and $[6,9,16,30]$, respectively.
By comparing the results with different $r_i$s, we can safely draw a conclusion that the
proposed method is robust to the choices of the $r_i$s.

 Our experiments also show that the results will suffer from  zigzag boundaries possibly if
all the radius of the radial functions are too large. On the other hand, the segmentation results will
have nonconvex boundary with small absolute curvature if all the radius of the radial functions are too small.
Therefore,  one only needs to use about four radial functions with radius between $4$ to $30$ to save computational cost for
real applications.

\section{Conclusion and future work}\label{sec:Con}
This this paper, we present a novel binary representation for convex shapes. It uses an inequality constraint on the indicator function.
This representation has two advantages. Firstly, It is a very general method which is independent of the dimension of the shape. Secondly,
the corresponding model with the proposed convexity constraint is very simple and easy to solve.

In the future, we will continue the research on this topic, such as
convexity representation methods, algorithms and applications.
Firstly, we will extend the proposed representation method for single convex object to the representation for multiple convex objects.
Experiments on 3D image data is on the way.

\section*{Appendix}
{\bf Proof of Theorem \ref{thm1}}~
    It is well known that $D$  is convex if and only if there is a
    hyper-plane $S(x)=\{y|N(x)\cdot y=N(x)\cdot x\}$ such that $D$ locates on one side of $S(x)$ for all $x\in\partial D$,
    where $N(x)$ is the normal vector of the hyper-plane.
    Without loss of generality, we assume
    \begin{equation}
      D\subset S^-(x)=\{y|N(x)\cdot y\leq N(x)\cdot x\}. \label{eq:ap1}
    \end{equation}
    It is obvious that $D\bigcap B_r(x)\subset S^-(x)\bigcap B_r(x)$ for $x\in\partial D$. Therefore, $D$ is
    convex if and only if
    \begin{align}
       \mathcal{A}(D\bigcap B_r(x))&\leq\mathcal{A}(B_r(x)\bigcap S^-(x))\notag\\
       &= \frac{1}{2}\mathcal{A}(B_r(x)), x~\text{on}~\partial D,
    \end{align}
    which is equivalent to $\mathcal{A}(D^c\bigcap B_r(x))\geq \frac{1}{2}\mathcal{A}(B_r(x))$.

{\bf Proof of Corollary \ref{cor1}}  We firstly prove the relation in (\ref{eq:char4}). According to the assumptions of $b_r$,
it is obvious that
\begin{equation}
\int_{B_r(x)\bigcap S^+(x)} b_r(y)dy=\frac{1}{2}\int_{B_r(x)} b_r(y)dy=\frac{1}{2}
\end{equation}
for $x \in \partial D$,
where $S^+(x)=\{y|N(x)\cdot y\geq N(x)\cdot x\}$.
Thus we can obtain easily
\begin{equation}
    \int_{B_r(x)\bigcap D^c}b_r(y-x)dy\geq \int_{B_r(x)\bigcap S^+} b_r(y-x)dy =\frac{1}{2}.\notag
\end{equation}
For $x\in D^c$, the relation (\ref{eq:char4}) is true obviously.

Therefore, the inequality in (\ref{eq:char5}) is true for $x\in D$ since $u(x)=1$ by (\ref{eq:char4}).
For $x\in D^c$, the inequality is also true obviously because $u(x)=0$.
{\small
\bibliographystyle{plain}
\bibliography{egbib}

\begin{thebibliography}{10}

\bibitem{Bae2017Augmented}
Egil Bae, Xue-Cheng Tai, and Wei Zhu.
\newblock Augmented {L}agrangian method for an {E}uler's elastica based
  segmentation model that promotes convex contours.
\newblock {\em Inverse Problems and Imaging}, 11(1):1--23, 2017.

\bibitem{Chan2005Level}
Tony~F. Chan and Wei Zhu.
\newblock Level set based shape prior segmentation.
\newblock In {\em IEEE Conference on Computer Vision and Pattern Recognition},
  volume~2, pages 1164--1170, 2005.

\bibitem{Cremers2003Towards}
Daniel Cremers and Nir Sochen.
\newblock Towards recognition-based variational segmentation using shape priors
  and dynamic labeling.
\newblock In {\em International Conference on Scale Space Methods in Computer
  Vision}, pages 388--400, 2003.

\bibitem{elsey2009analogue}
Matthew Elsey and Selim Esedoglu.
\newblock Analogue of the total variation denoising model in the context of
  geometry processing.
\newblock {\em Multiscale Modeling \& Simulation}, 7(4):1549--1573, 2009.

\bibitem{esedog2015threshold}
Selim Esedo\-g Lu and Felix Otto.
\newblock Threshold dynamics for networks with arbitrary surface tensions.
\newblock {\em Communications on pure and applied mathematics}, 68(5):808--864,
  2015.

\bibitem{gorelick2017multi}
Lena Gorelick and Olga Veksler.
\newblock Multi-object convexity shape prior for segmentation.
\newblock In {\em International Workshop on Energy Minimization Methods in
  Computer Vision and Pattern Recognition}, pages 455--468. Springer, 2017.

\bibitem{Gorelick2017Convexity}
Lena Gorelick, Olga Veksler, Yuri Boykov, and Claudia Nieuwenhuis.
\newblock Convexity shape prior for binary segmentation.
\newblock {\em IEEE Transactions on Pattern Analysis and Machine Intelligence},
  39(2):258--270, 2017.

\bibitem{gulshan2010geodesic}
Varun Gulshan, Carsten Rother, Antonio Criminisi, Andrew Blake, and Andrew
  Zisserman.
\newblock Geodesic star convexity for interactive image segmentation.
\newblock In {\em IEEE Conference on Computer Vision and Pattern Recognition},
  pages 3129--3136, 2010.

\bibitem{isack2018k}
Hossam Isack, Lena Gorelick, Karin Ng, Olga Veksler, and Yuri Boykov.
\newblock K-convexity shape priors for segmentation.
\newblock In {\em European Conference on Computer Vision}, pages 36--51, 2018.

\bibitem{isack2016hedgehog}
Hossam Isack, Olga Veksler, Milan Sonka, and Yuri Boykov.
\newblock Hedgehog shape priors for multi-object segmentation.
\newblock In {\em IEEE Conference on Computer Vision and Pattern Recognition},
  pages 2434--2442, 2016.

\bibitem{Leventon2003Statistical}
Michael~E. Leventon, W.~Eric~L. Grimson, and Olivier Faugeras.
\newblock Statistical shape influence in geodesic active contours.
\newblock In {\em IEEE Embs International Summer School on Biomedical Imaging},
  pages 316--322, 2000.

\bibitem{li2019variational}
Lingfeng Li, Shousheng Luo, Xue-Cheng Tai, and Jiang Yang.
\newblock A variational convex hull algorithm.
\newblock In {\em International Conference on Scale Space and Variational
  Methods in Computer Vision}, pages 224--235. Springer, 2019.

\bibitem{liu2011fast}
Jun Liu, Xue-cheng Tai, Haiyang Huang, and Zhongdan Huan.
\newblock A fast segmentation method based on constraint optimization and its
  applications: Intensity inhomogeneity and texture segmentation.
\newblock {\em Pattern Recognition}, 44(9):2093--2108, 2011.

\bibitem{Luo:2018tl}
Shousheng Luo and Xue-Cheng Tai.
\newblock Convex shape priors for level set representation.
\newblock {\em arXiv preprint arXiv:1811.04715}, 2018.

\bibitem{luo2019Cnvx}
Shousheng Luo, Xue-Cheng Tai, Limei Huo, Yang Wang, and Roland Glowinsiki.
\newblock Convex shape prior for multi-object segmentation using a single level
  set function.
\newblock In {\em International Conference on Computer Vision}, pages 613--621,
  2019.

\bibitem{Mahabadi2015Segment}
Rabeeh~Karimi Mahabadi, Christian Hane, and Marc Pollefeys.
\newblock Segment based {3D} object shape priors.
\newblock In {\em IEEE Conference on Computer Vision and Pattern Recognition},
  2015.

\bibitem{mirpalpar2007}
Michele Miranda, Diego Pallara, Fabio Paronetto, and Marc Preunkert.
\newblock Short-time heat flow and functions of bounded variation in
  {$\mathbb{R}^N$}.
\newblock {\em Annales de la Facult\'e des Sciences de Toulouse:
  Math\'ematiques}, 16(1):125--145, 2007.

\bibitem{mirikharaji2018star}
Zahra Mirikharaji and Ghassan Hamarneh.
\newblock Star shape prior in fully convolutional networks for skin lesion
  segmentation.
\newblock In {\em International Conference on Medical Image Computing and
  Computer-Assisted Intervention}, pages 737--745. Springer, 2018.

\bibitem{Rother:2004}
Carsten Rother, Vladimir Kolmogorov, and Andrew Blake.
\newblock "grabcut": Interactive foreground extraction using iterated graph
  cuts.
\newblock {\em ACM Transactions on Graphics (TOG)}, 23(3):309--314, 2004.

\bibitem{royer2016convexity}
Loic~A. Royer, David~L. Richmond, Carsten Rother, Bjoern Andres, and Dagmar
  Kainmueller.
\newblock Convexity shape constraints for image segmentation.
\newblock In {\em IEEE Conference on Computer Vision and Pattern Recognition},
  pages 402--410, 2016.

\bibitem{toranzos2004sets}
Fausto~A Toranzos and Ana~Forte Cunto.
\newblock Sets expressible as finite unions of star shaped sets.
\newblock {\em Journal of Geometry}, 79(1-2):190--195, 2004.

\bibitem{Ukwatta2013Efficient}
Eranga Ukwatta, Jing Yuan, Wu~Qiu, Martin Rajchl, and Aaron Fenster.
\newblock Efficient convex optimization-based curvature dependent contour
  evolution approach for medical image segmentation.
\newblock In Sebastien Ourselin and David~R Haynor, editors, {\em Medical
  Imaging 2013: Image Processing}, volume 8669, pages 866--902, 2013.

\bibitem{veksler2008star}
Olga Veksler.
\newblock Star shape prior for graph-cut image segmentation.
\newblock In {\em European Conference on Computer Vision}, pages 454--467.
  Springer, 2008.

\bibitem{vicente2008graph}
Sara Vicente, Vladimir Kolmogorov, and Carsten Rother.
\newblock Graph cut based image segmentation with connectivity priors.
\newblock In {\em 2008 IEEE Conference on Computer Vision and Pattern
  Recognition}, pages 1--8. IEEE, 2008.

\bibitem{wang2017efficient}
Dong Wang, Haohan Li, Xiaoyu Wei, and Xiao-Ping Wang.
\newblock An efficient iterative thresholding method for image segmentation.
\newblock {\em Journal of Computational Physics}, 350:657--667, 2017.

\bibitem{wang2009edge}
Jie Wang, Lili Ju, and Xiaoqiang Wang.
\newblock An edge-weighted centroidal {Voronoi} tessellation model for image
  segmentation.
\newblock {\em IEEE Transactions on Image Processing}, 18(8):1844--1858, 2009.

\bibitem{Yan2018}
Shi Yan, Xue-Cheng Tai, Jun Liu, and Haiyang Huang.
\newblock Convexity shape prior for level set based image segmentation method.
\newblock {\em arXiv preprint arXiv:1805.08676}, 2018.

\bibitem{yang2017level}
Cong Yang, Xue Shi, Donglan Yao, and Chunming Li.
\newblock A level set method for convexity preserving segmentation of cardiac
  left ventricle.
\newblock In {\em International Conference on Image Processing}, pages
  2159--2163, 2017.

\bibitem{yuan2012efficient}
Jing Yuan, Wu~Qiu, Eranga Ukwatta, Martin Rajchl, Yue Sun, and Aaron Fenster.
\newblock An efficient convex optimization approach to 3d prostate mri
  segmentation with generic star shape prior.
\newblock {\em Prostate MR Image Segmentation Challenge, MICCAI}, 7512:82--89,
  2012.

\bibitem{yuan2012fast}
Jing Yuan, Eranga Ukwatta, Xue-Cheng Tai, A~Fenster, and C~Schnoerr.
\newblock A fast global optimization-based approach to evolving contours with
  generic shape prior.
\newblock Technical report, UCLA, 2012.

\bibitem{Zafari2016}
Sahar Zafari, Tuomas Eerola, Jouni Sampo, Heikki K$\ddot{a}$lvi$\ddot{a}$inen,
  and Heikki Haario.
\newblock Segmentation of partially overlapping convex objects using branch and
  bound algorithm.
\newblock In {\em Asian Conference on Computer Vision}, pages 76--90, 2016.

\end{thebibliography}
}
\end{document}